\documentclass[11pt]{article}
\usepackage{amsmath,amsthm,amsfonts}
\usepackage[dvips]{epsfig}
\usepackage{fullpage}
\usepackage{float}

\newtheorem{theorem}{Theorem}

\def\mod#1 #2{#1\ ({\rm mod}\ #2)}

\title{Hamming Distance for Conjugates}

\author{Jeffrey Shallit\\
School of Computer Science\\
University of Waterloo\\
Waterloo, Ontario  N2L 3G1  Canada\\
{\tt shallit@graceland.uwaterloo.ca}}

\begin{document}

\maketitle

\begin{abstract}
     Let $x, y$ be strings of equal length.  
The {\it Hamming distance} $h(x,y)$
between $x$ and $y$ is the number of positions in which $x$ and $y$
differ.  If $x$ is a cyclic shift of $y$,
we say $x$ and $y$ are {\it conjugates}.
We consider $f(x,y)$,
the Hamming distance between the conjugates $xy$ and $yx$.
Over a binary alphabet $f(x,y)$ is always even, and must satisfy a further
technical condition.  By contrast, over an
alphabet of size $3$ or greater, $f(x,y)$ can take any
value between $0$ and $|x|+|y|$, except $1$; furthermore, we can
always assume that the smaller string has only one type of letter. 
\end{abstract}

%
%

\section{Introduction}
     Let $x, y$ be strings of equal length.  We define the
{\it Hamming distance} $h(x,y)$
between $x$ and $y$ to be the number of positions in which $x$ and $y$
differ \cite{Hamming:1950}.
Thus, for example, $h({\tt seven}, {\tt three}) = 4$.  
If $x$ is a cyclic shift of $y$, we say $x$ and $y$ are {\it conjugates}.
Alternatively, $x$ is a conjugate of $y$ if there exist strings $u,v$
with $x = uv$ and $y = vu$.  For example, $x = {\tt enlist}$ and
$y = {\tt listen}$ are conjugates; take $u = {\tt en}$, $v = {\tt list}$.

     In this paper we consider the Hamming distance for conjugates.
For all strings $x, y$ (not necessarily of the same length) we define
$f(x,y) = h(xy, yx)$.
If $f(x,y) = 0$, then $xy = yx$, and this equation holds if and only
if both $x$ and $y$ are powers of some string $z$
\cite{Lyndon&Schutzenberger:1962}.  On the other hand,
there are no solutions to the equation $f(x,y) = 1$.  For as a referee
of an earlier version of this paper noted,
if $h(xy, yx) = 1$, then we can write $xy = uav$, $yx = ubv$ for some
strings $u, v$ and letters $a, b$ with $a \not= b$.  Thus $xy$ contains
one more $a$ than $yx$, which is impossible, since $xy$ and $yx$ clearly
contain the same number of occurrences of each letter.

      These two examples suggest trying to determine for which $k$ the
equation $f(x,y) = k$ is solvable.  Over an
alphabet of size $3$ or greater, we show that $f(x,y)$ can take any
value between $0$ and $|x|+|y|$, except $1$; furthermore, we can
always assume that the smaller string has only one type of letter. 
However, over a binary alphabet, $f(x,y)$ is always even, and must
also satisfy a further technical condition.

\section{Hamming distance for non-binary alphabets}

\begin{theorem}
Let $\Sigma$ be an alphabet with at least $3$ letters, say
${\tt 0}$, ${\tt 1}$, ${\tt 2}$. 
Suppose $m, n, k$ are integers with $1 \leq m \leq n$, 
and $0 \leq k \leq m+n$.
\begin{itemize}

\item[(a)]
If $m < n$,
there exist strings $x, y \in \Sigma^*$ with $|x| = m$, $|y| = n$,
such that $h(xy,yx) = k$ if and only if $k \not= 1$.

\item[(b)] If
$m = n$, then there exist strings $x,y \in \Sigma^*$ with $|x| = m$, $|y| = n$,
such that $h(xy,yx) = k$ if and only if $k$ is even.
\end{itemize}

Furthermore, in both cases, we can always choose $x = {\tt 0}^m$.
\end{theorem}

\begin{proof}

We define $x = {\tt 0}^m$ and $y = s(m,n,k)$, where $s(m,n,k)$ is defined
by the following recursion:

\begin{eqnarray}
s(m,n,2t) &=& {\tt 0}^{n-t} {\tt 1}^t,
\quad\quad \text{if } 0 \leq t \leq m \leq n \label{eq1} \\
s(m,n,2t+1) &=& {\tt 0}^{n-m-1} {\tt 1 0}^{m-t} {\tt 1}^{t-1} {\tt 2},
\quad\quad \text{if }
	1 \leq t \leq m < n 
	\label{eq2} \\
s(m,n,k) &=& {\tt 0}^{n+m-k} s(m,k-m,k),
\quad\quad \text{if } 2m \leq k \leq m+n 
	\label{eq3} \\
s(m,n,m+n) &=& \begin{cases}
	({\tt 1}^m {\tt 2}^m)^j {\tt 1}^r,
	& \text{ if } n = 2mj+r,\ 0 \leq r \leq m \\
	({\tt 1}^m {\tt 2}^m)^j {\tt 1}^m {\tt 2}^r, & \text{ if } n = (2m+1)j + r,\ 0 \leq r \leq m \\
	\end{cases}
	\label{eq4} 
\end{eqnarray}

     First we prove that that these identities suffice to calculate
$s(m,n,k)$ for $0 \leq k \leq m+n$, $k \not=1 $, when $m < n$ and
for $0 \leq k \leq m+n$, $k$ even, when $m = n$.

     Suppose $m < n$.  Then there are two cases:  either $k \leq 2m+1$ or
$k > 2m+1$.  Suppose $k \leq 2m+1$.    If $k$ is even, in which case $k = 2t$,
$0 \leq t \leq m$, we use Eq.~(\ref{eq1}).  If
$k$ is odd, in which case $k = 2t+1$, $1 \leq t \leq m$, we
use Eq.~(\ref{eq2}).  

    Now suppose $k > 2m+1$.  Then if $k < m + n$, we use Eq.~(\ref{eq3}),
which reduces the case to one where $k = m+n$.  In this latter case,
we use Eq.~(\ref{eq4}).  

     If $m = n$, then we use Eq.~(\ref{eq1}) if $k \leq 2m$, and
Eqs.~(\ref{eq3}) and (\ref{eq4}) if $k > 2m$.    Thus the identities
(\ref{eq1})--(\ref{eq4}) cover all the cases.   
Furthermore, if $|x| = |y|$ then
$h(xy,yx) = 2 h(x,y)$, so $k$ must be even.

    Now we prove that $f({\tt 0}^m, s(m,n,k)) = k$.   We start with (\ref{eq1}).
Comparing ${\tt 0}^m {\tt 0}^{n-k} {\tt 1}^k$ with
${\tt 0}^{n-k} {\tt 1}^k {\tt 0}^m$, we see that since
$m \geq k$, each $\tt 1$ is paired with a $\tt 0$ in the other string, and all
other symbols are $\tt 0$, so $f({\tt 0}^m , {\tt 0}^{n-k} {\tt 1}^k) = 2k$.

    Now consider (\ref{eq2}).  By comparing 
$xy = {\tt 0}^m {\tt 0}^{n-m-1} {\tt 1} {\tt 0}^{m-k} {\tt 1}^{k-1} {\tt 2}$
with
$yx = {\tt 0}^{n-m-1} {\tt 1 0}^{m-k} {\tt 1}^{k-1} {\tt 2 0}^m$ we see that,
since
$k \leq m$, the last $k$ symbols of $xy$ are different from $\tt 0$, while
the last $k$ symbols of $yx$ are $\tt 0$.
The block ${\tt 1}^{k-1} {\tt 2}$ in $yx$
is matched against ${\tt 0}^{k-1} {\tt 1}$ in $xy$.
And the first $\tt 1$ in $yx$ is
matched against $\tt 0$ in $xy$.  The total is $2k+1$ mismatches, as needed.

    Now consider (\ref{eq3}).   Adding $0$'s here to the front of 
$y = s(m, k-m,k)$ does not change the number of mismatches.

    Finally, consider (\ref{eq4}).  In this case, the $\tt 0$'s in $xy$
match against $\tt 1$'s in $yx$.  The alternating blocks of $\tt 1$'s and 
$\tt 2$'s in
$xy$ either match against blocks of the other symbol in $yx$
($\tt 1$'s against $\tt 2$'s and $\tt 2$'s against $\tt 1$'s),
or against the $\tt 0$'s
at the end of $yx$.  Thus every symbol mismatches, and there are $m+n$
of them.

     This completes the proof.
\end{proof}

\smallskip

\noindent{\bf Examples.}  We consider some examples.  

\begin{eqnarray*}
s(10, 45, 55) &=& {\tt 1}^{10} {\tt 2}^{10} {\tt 1}^{10} {\tt 2}^{10} {\tt 1}^5 ; \\
s(10, 11, 21) &=& {\tt 1}^{10} {\tt 2} ; \\
s(10, 10, 20) &=& {\tt 1}^{10} ; \\
s(10, 20, 12) &=& {\tt 0}^{14} {\tt 1}^6 ;\\
s(10, 20, 13) &=& {\tt 0}^{9} {\tt 1 0}^4 {\tt 1}^5 {\tt 2} .\\
\end{eqnarray*}

\section{The case of a binary alphabet}

      Suppose $x, y$ are strings over the alphabet $\lbrace {\tt 0, 1} \rbrace$.
We first observe that $f(x,y)$ must always be even.  For, as the
referee of an earlier version observed, if $xy$ and $yx$ differ at an
odd number of positions, then one of them must contain at least one more
${\tt 0}$ than the other.  

However, to completely characterize the solutions to $f(x,y) = k$,
there is one additional condition that needs to be imposed:

\begin{theorem}
     Let $m, n$ be integers with $1 \leq m \leq n $.  Then
there exist binary strings $x, y$ with $f(x,y) = k$ if and only
if
\begin{itemize}
\item[(a)] $0 \leq k \leq m+n$;

\item[(b)] $k$ is even;

\item[(c)] $k \leq m + n - \gcd(m,n)$ if 
$(m+n)/\gcd(m,n)$ is odd.
\end{itemize}
\end{theorem}

\begin{proof}
Suppose $f(x,y) = k$ is solvable.  We have already seen that conditions
(a) and (b) must hold.  By comparing $xy$ to $yx$ we
see that each symbol is potentially related to $(m+n)/\gcd(m,n) - 1$
other symbols.  
For example, writing $xy = z$, and using indexing beginning at $0$,
we see that $z[0]$ is the first symbol
of $xy$ and $z[m]$ is the first symbol of $yx$.  The $m$'th symbol
of $xy$ is equal to the $2m \bmod (m+n)$'th symbol of $yx$, and so
forth.  That is, the positions of $x$ and $y$ split into
$\gcd(m,n)$ cycles of length $(m+n)/\gcd(m,n)$; adjacent elements of
a cycle line up with each other in $xy$ and $yx$.  

If a cycle is of even length, then over a binary alphabet we can force
all the symbols to disagree, by choosing them to be $\tt 0$ and $\tt 1$ 
alternately.  If a cycle is of odd length, this is impossible.  
More precisely, the number of adjacent pairs that differ in a cycle
must be even.  

Therefore, if $(m+n)/\gcd(m,n)$ is odd, at most $(m+n)/\gcd(m,n) - 1$
pairs of any cycle can disagree.  Since there are $\gcd(m,n)$ cycles,
the highest Hamming weight we can achieve is $m+n - \gcd(m,n)$.
Thus conditions (a)--(c) must hold.

Now suppose conditions (a)--(c) hold.  We show how to construct
$x$, $y$ such that $f(x,y) = k$.  
Define $g$ to be $(m+n)/\gcd(m,n)$ if this quantity is odd; otherwise
let $g = (m+n)/\gcd(m,n) - 1$.
Using the division theorem, divide $k$ by $g$, obtaining a quotient
$q$ and a remainder $r$.  Since $k$ and $g$ are even, so is $r$.
In the first $q$ of the cycles, let the symbols alternate between $\tt 0$
and $\tt 1$.
In the $(q+1)$'th cycle, let the first $r$ symbols alternate and set
all the remaining symbols to be the same.  The resulting string, now 
split up between $x$ and $y$, now has $qg+r = k$ positions where $xy$
fails to match $yx$, as desired.
\end{proof}

\bigskip

\noindent{\bf Example.}   Suppose $(m,n,k) = (6,9,10)$.  In this
case, the cycles are 
\begin{eqnarray*}
c_0 &=& (z[0],z[6],z[12],z[3],z[9]); \\
c_1 &=& (z[1],z[7],z[13],z[4],z[10]); \\
c_2 &=& (z[2],z[8],z[14],z[5],z[11]). \\
\end{eqnarray*}
Now $(m+n)/\gcd(m,n) = 5$, which is odd.  Hence each cycle can give
us at most $4$ mismatches.    We generate $4$ mismatches with each
of the first two cycles by alternating $0$ and $1$,
and generate $2$ mismatches with the last cycle.
This gives
\begin{eqnarray*}
c_0 &=& (z[0],z[6],z[12],z[3],z[9]) = {\tt  (0,1,0,1,0)}; \\
c_1 &=& (z[1],z[7],z[13],z[4],z[10]) = {\tt (0,1,0,1,0)};  \\
c_2 &=& (z[2],z[8],z[14],z[5],z[11]) = {\tt (0,1,0,0,0)}.\\
\end{eqnarray*}
This gives $z = {\tt 000111000110000}$ and so
$x = {\tt 000110}$ and $y= {\tt 111000000}$.

\bigskip

It is interesting to note that, in contrast to the case of large alphabets,
in the binary case, 
even if $x$ and $y$ exist with $|x| = m$, $|y| = n$, $m \leq n$,
and  $h(x,y) = k$, it may not be possible
to achieve this by choosing $x = {\tt 0}^m$.  For example, for
$(m,n,k) = (3,5,8)$, the lexicographically smallest solution is
$x = {\tt 010}$, $y = {\tt 10101}$.

\section{Acknowledgments}

I am very grateful to the anonymous referee for substantially simplifying
the proofs of some of the claims.

\end{document}